\documentclass[11pt]{article}
\usepackage{amsmath,latexsym,epsf,mathabx}
\usepackage{epsfig,graphicx,picins,picinpar,subfigure}
\usepackage{pstricks}
\usepackage{fancyvrb}
\usepackage[all]{xy}

\usepackage{ulem}
\usepackage{mathrsfs}
\usepackage{soul}
\usepackage{color}
\begin{document}

\newcommand{\nc}{\newcommand}
\def\PP#1#2#3{{\mathrm{Pres}}^{#1}_{#2}{#3}\setcounter{equation}{0}}
\def\ns{$n$-star}\setcounter{equation}{0}
\def\nt{$n$-tilting}\setcounter{equation}{0}
\def\Ht#1#2#3{{{\mathrm{Hom}}_{#1}({#2},{#3})}\setcounter{equation}{0}}
\def\qp#1{{${(#1)}$-quasi-projective}\setcounter{equation}{0}}
\def\mr#1{{{\mathrm{#1}}}\setcounter{equation}{0}}
\def\mc#1{{{\mathcal{#1}}}\setcounter{equation}{0}}
\def\HD{\mr{Hom}_{\mc{D}}}
\def\HC{\mr{Hom}_{\mc{C}}}
\def\AdT{\mr{Add}_{\mc{T}}}
\def\adT{\mr{add}_{\mc{T}}}
\def\Kb{\mc{K}^b(\mr{Proj}R)}
\def\kb{\mc{K}^b(\mc{P}_R)}
\def\AdpC{\mr{Adp}_{\mc{C}}}

\newtheorem{Th}{Theorem}[section]
\newtheorem{Def}[Th]{Definition}
\newtheorem{Lem}[Th]{Lemma}
\newtheorem{Pro}[Th]{Proposition}
\newtheorem{Cor}[Th]{Corollary}
\newtheorem{Rem}[Th]{Remark}
\newtheorem{Not}[Th]{Notation}
\newtheorem{Sc}[Th]{}
\def\Pf#1{{\noindent\bf Proof}.\setcounter{equation}{0}}
\def\>#1{{ $\Rightarrow$ }\setcounter{equation}{0}}
\def\<>#1{{ $\Leftrightarrow$ }\setcounter{equation}{0}}
\def\bskip#1{{ \vskip 20pt }\setcounter{equation}{0}}
\def\sskip#1{{ \vskip 5pt }\setcounter{equation}{0}}
\def\mskip#1{{ \vskip 10pt }\setcounter{equation}{0}}
\def\bg#1{\begin{#1}\setcounter{equation}{0}}
\def\ed#1{\end{#1}\setcounter{equation}{0}}
\def\KET{T^{^F\bot}\setcounter{equation}{0}}
\def\KEC{C^{\bot}\setcounter{equation}{0}}

\renewcommand{\thefootnote}{\fnsymbol{footnote}}
\setcounter{footnote}{0}
%
%


\title{\bf Star modules with respect to balanced pair
\thanks{Supported by the National Natural Science Foundation of China (Grants No.11801004), the Startup Foundation for Introducing Talent of AHPU (Grant No.2017YQQ016)
and the Top talent project of AHPU in 2020. }
}
\smallskip
\author{ Peiyu Zhang and  Dajun Liu\\ 
\footnotesize ~E-mail:~zhangpy@ahpu.edu.cn ~ldjnnu2017004@163.com;\\
\footnotesize School of Mathematics and Physics,  Anhui Polytechnic University, Wuhu, China\\
Jiaqun Wei\\ 
\footnotesize ~E-mail:~ weijiaqun@njnu.edu.cn; \\
\footnotesize School of Mathematics Science, Nanjing Normal University, Nanjing, China
}
\date{}
\maketitle
\baselineskip 15pt
%
%
\begin{abstract}
\vskip 10pt%

In this article, firstly, we introduce the notion of star modules with respect to a balanced pair and obtain some properties.
We mainly give the relationship between $n$-${\mathscr X}$ star modules and $n$-${\mathscr X}$ tilting modules \cite{zpy},
and a new characterization of $n$-${\mathscr X}$ tilting modules.

\mskip\




\end{abstract}
%
\vskip 30pt
\section{Introduction and Preliminaries}

Let $R$ be an ring, denote all left $R$-modules by $\mathrm{Mod}$$R$. Throughout this paper, all $R$-modules are left $R$-module.
Assume that $\mathscr C$ is an subcategory of $\mathrm{Mod}$$R$. For any left $R$-module $M$, if there is a $\mathscr C$-exact complex
$C^{\bullet}=:$ $\cdots \longrightarrow C_{2}\longrightarrow C_{1}\longrightarrow C_{0}\longrightarrow M$
with $C_{i}\in \mathscr C$ for $i\geq0$ (i.e., it is exact by applying the functor $\mathrm{Hom}_{R}(C,-)$ for each $C\in\mathscr C$),
then we call the complex being $\mathscr C$-resolution of $M$.
Recall that an subcategory $\mathscr B$ of $\mathrm{Mod}$$R$ is said to be contravariantly
finite, if for any left $R$-module $A$, it has a right $\mathscr B$-approximation \cite{AR}, i.e., there is
a homomorphism $f$ : $B\to A$ for some $B\in \mathscr B$ such
that $\mathrm{Hom}_{R}(B',f)$ is surjective for any $B'\in \mathscr B$.
Dually, we have the definition of the right $\mathscr C$-resolution and covariantly finite subcategory.

The notion of the balanced pair was first introduced by Chen \cite{CHEN}, the author mainly given a triangle-equivalence
between two homotopy categories of complexes, see \cite[Theorem~4.1]{CHEN}.  Furthermore, the author proved that there exists a
triangle-equivalence between the homotopy category of Gorenstein projective modules and the
homotopy category of Gorenstein injective modules when $R$ is a left-Gorenstein ring.

\bg{Def}$\mathrm{\cite{CHEN}}$\label{BP}%
A pair $({\mathscr X}$, ${\mathscr Y})$ of additive subcategories in an Abelian category ${\mathscr A}$ is called a balanced pair if the following
conditions are satisfied:

$(1)$ the subcategory ${\mathscr X}$ is contravariantly finite and ${\mathscr Y}$ is covariantly finite;

$(2)$ for each object M, there is an ${\mathscr X}$-resolution $X^{\bullet}\longrightarrow M$ such that it is exact by applying the
functors $\mathrm{Hom}_{{\mathscr A}}(-,~Y )$ for all $Y\in{\mathscr Y}$;

$(3)$ for each object N, there is a ${\mathscr Y}$-coresolution $N\longrightarrow Y^{\bullet}$ such that it is exact by applying the
functors $\mathrm{Hom}_{{\mathscr A}}(X,~-)$ for all $X\in{\mathscr X}$.
\ed{Def}

Let $\mathscr A$ be an Abelian category. A contravariantly finite subcategory $\mathscr C$ of $\mathscr A$ is admissible provided that each right
$\mathscr C$-approximation is surjective. It is equivalent to that any ${\mathscr X}$-exact complex is indeed exact.
Dually, the notion of coadmissible is defined.
Recalled that a balanced pair $({\mathscr X}$, ${\mathscr Y})$ is admissible, if ${\mathscr X}$ is admissible
(or,${\mathscr Y}$ is coadmissible, see the corollary 2.3 in \cite{CHEN}). On the research of balanced pairs,
for example,  Zheng constructed a recollement of triangulated categories through a complete balanced pair, see \cite[Theorem~2.12]{zyf}.
Li \cite{LI} introduced the notion of cotorsion pair \cite{Enochs,H,HI,KS} with respect to a balanced pair, and proved that the right ${\mathscr X}$-singularity category coincides with
the left ${\mathscr Y}$-singularity category if the
${\mathscr X}$-resolution dimension of ${\mathscr Y}$ and the ${\mathscr Y}$-coresolution dimension of ${\mathscr X}$ are finite.

In this paper, we also assume that $({\mathscr X}$, ${\mathscr Y})$ is an admissible balanced pair in $\mathrm{Mod}$$R$.
Whenever the subcategory ${\mathscr X}$ (respectively, ${\mathscr Y}$) appears, we always suppose
that it is a part of the admissible balanced pair $({\mathscr X}$, ${\mathscr Y})$.
Note that the functor $\mathrm{Ext}^{i}_{{\mathscr A}}(-,~-)$ is based on the classical balanced pair $(\mathrm{Proj} R,~\mathrm{Inj} R)$, it induces
an isomorphism of cohomology groups whether we take a projective resolution of the first variable or take
an injective coresolution of the second variable. From this viewpoint, we have the functor $\mathrm{Ext}^{i}_{{\mathscr X}}(-,~-)$ with respect
to admissible balanced pair $({\mathscr X}$, ${\mathscr Y})$. By Lemma 2.1 \cite{CHEN}, the functor is well-define.

For an $R$-module $Q$, we have some notions which are widely used in the relative tilting theory.

$\mathrm{Add}Q$ =: $\{K|$ there exists a $R$-module $J$ such that $K\bigoplus J\cong Q^{(I)}$ for some set $I. \}$

$\mathrm{Pres}^{n}_{\mathscr X}Q$ =: $\{L|$ there is $\mathscr X$-exact sequence $Q_{n}\longrightarrow \cdots\longrightarrow Q_{1}\longrightarrow L\longrightarrow0$
with $Q_{i}\in\mathrm{Add}Q$ for all $1\leq i\leq n.\}$

$Q^{\mathscr X \perp_{i\geq1}}$=: $\{M|$ $\mathrm{Ext}^{i\geq1}_{\mathscr X}(Q,~M)=0.\}$

\section{Relative star modules}

In this section, we mainly give the definition of $n$-${\mathscr X}$ module and some properties.
Firstly, we introduce the notion of $n$-${\mathscr X}$ quasi-projective module.

\bg{Def}\label{x-quasi}
An R-module Q is called $n$-${\mathscr X}$ quasi-projective if for any ${\mathscr X}$-exact sequence
$0\longrightarrow A\longrightarrow Q_{0}\longrightarrow B\longrightarrow0$ with $Q_{0}\in \mathrm{Add}Q$ and $A\in \mathrm{Pres}^{n-1}_{\mathscr X}Q$,
the sequence is Q-exact.
\ed{Def}

If it satisfies condition $\mathrm{Pres}^{n-1}_{\mathscr X}Q=\mathrm{Pres}^{n}_{\mathscr X}Q$ for $R$-module $Q$, then we can obtain a
characterization of $n$-${\mathscr X}$ quasi-projective modules as follows.

\bg{Pro}\label{char}
For an R-module Q, if $\mathrm{Pres}^{n-1}_{\mathscr X}Q=\mathrm{Pres}^{n}_{\mathscr X}Q$, then the following statements are equivalent.

$(1)$ Q is an $n$-${\mathscr X}$ quasi-projective module;

$(2)$ If $\xymatrix{0\ar[r]&A\ar[r]&Q_{0}\ar[r]&B\ar[r]&0}$ is $\mathscr X$-exact with $Q_{0}\in \mathrm{Add}Q$
and $B\in \mathrm{Pres}^{n}_{\mathscr X}Q$, then $A\in \mathrm{Pres}^{n-1}_{\mathscr X}Q$ if and only if the sequence is $Q$-exact.
\ed{Pro}

\Pf. By the definition of $n$-${\mathscr X}$ quasi-projective module, we only need to verify that $(1)\Rightarrow(2)$.

If $A\in \mathrm{Pres}^{n-1}_{\mathscr X}Q$, then the result is clear. On the other hand, let the above sequence be $Q$-exact.
Since $B\in \mathrm{Pres}^{n}_{\mathscr X}Q$, there is an ${\mathscr X}$-exact sequence
$\xymatrix{0\ar[r]&B_{1}\ar[r]&Q_{B}\ar[r]&B\ar[r]&0}$ with $Q_{B}\in \mathrm{Add}Q$ and $B_{1}\in \mathrm{Pres}^{n-1}_{\mathscr X}Q$.
Obviously, it is $Q$-exact by (1). Thus we have that $A\bigoplus Q_{B}\cong B_{1}\bigoplus Q_{0} \in\mathrm{Pres}^{n-1}_{\mathscr X}Q
=\mathrm{Pres}^{n}_{\mathscr X}Q$ by the lemma 2.3 in \cite{WEI2}. There is an ${\mathscr X}$-exact sequence
$0\longrightarrow A'\longrightarrow Q'_{0}\longrightarrow A\bigoplus Q_{B}\longrightarrow0$ with $Q'_{0}\in \mathrm{Add}Q$ and
$A'\in \mathrm{Pres}^{n-1}_{\mathscr X}Q$. So we can obtain the following commutative diagram.

$$\xymatrix{
&0\ar[d] & 0\ar[d]&\\
&A'\ar@^{=}[r]\ar[d] & A'\ar[d]\\
0\ar[r] &B_{2}\ar[r]\ar[d] &Q'_{0}\ar[r]\ar[d] &A\ar@^{=}[d]\ar[r] &0\\
0\ar[r]&Q_{B}\ar[r]\ar[d] &A\bigoplus Q_{B}\ar[d]\ar[r]&A\ar[r] &0\\
&0&0&&}$$
It is easy to prove that all columns and rows in above diagram are $Q$-exact and ${\mathscr X}$-exact. Since $A'\in \mathrm{Pres}^{n-1}_{\mathscr X}Q$,
there is an ${\mathscr X}$-exact sequence $\xymatrix{0\ar[r]&A''\ar[r]&Q_{A'}\ar[r]&A'\ar[r]&0}$ with $Q_{A'}\in \mathrm{Add}Q$
and $A''\in \mathrm{Pres}^{n-2}_{\mathscr X}Q$. Consequently, we have the following commutative diagram.
$$\xymatrix{
&0\ar[d] & 0\ar[d]&\\
&A''\ar@^{=}[r]\ar[d] &A''\ar[d]\\
0\ar[r] &Q_{A'}\ar[r]^{(1~0)'}\ar[d]^{\pi_{1}} &Q_{A'}\bigoplus Q_{B}\ar[r]^{(0~1)}\ar[d]^{(i_{2}\pi_{1},~f)} &Q_{B}\ar@^{=}[d]\ar[r] &0\\
0\ar[r]&A'\ar[r]^{i_{2}}\ar[d] &B_{2}\ar[d]\ar[r]^{\pi_{2}}&Q_{B}\ar[r] &0\\
&0&0&&}$$
The existence of morphism $f$ is based on the fact that the third row in above diagram is $Q$-exact.
It follows from the second column that $B_{2}\in\mathrm{Pres}^{n-1}_{\mathscr X}Q$ since $A''\in \mathrm{Pres}^{n-2}_{\mathscr X}Q$.
From the second row in the first diagram, we have that $A\in\mathrm{Pres}^{n}_{\mathscr X}Q=\mathrm{Pres}^{n-1}_{\mathscr X}Q$.
\ \hfill $\Box$

\mskip\

Now, we give the definition of $n$-${\mathscr X}$ module and some useful properties.

\bg{Def}\label{xstar}
An R-module S is called an star module with respect to $({\mathscr X}$, ${\mathscr Y})$ (or, called $n$-$\mathscr X$ star module) if S is
$(n+1)$-$\mathscr X$ quasi-projective and $\mathrm{Pres}^{n}_{\mathscr X}S=\mathrm{Pres}^{n+1}_{\mathscr X}S$.
\ed{Def}

\bg{Pro}\label{123}
Let S be an $n$-$\mathscr X$ star module and the sequence $0\longrightarrow X\longrightarrow Y\longrightarrow Z\longrightarrow 0$
be $\mathscr X$-exact.

$(1)$ If both $X$ and $Y$ are in $\mathrm{Pres}^{n}_{\mathscr X}S$, then $Z$ is also in $\mathrm{Pres}^{n}_{\mathscr X}S$;

$(2)$ If the sequence is $S$-exact and $Z\in\mathrm{Pres}^{n}_{\mathscr X}S$, then $X\in\mathrm{Pres}^{n}_{\mathscr X}S$
if and only if $Y\in\mathrm{Pres}^{n}_{\mathscr X}S$;

$(3)$ If X, Y and Z are in $\mathrm{Pres}^{n}_{\mathscr X}S$, then the sequence is $S$-exact.

\ed{Pro}

\Pf. (1) Since $X, Y\in\mathrm{Pres}^{n}_{\mathscr X}S=\mathrm{Pres}^{n+1}_{\mathscr X}S$, there is two $\mathscr X$-exact sequences
$0\longrightarrow X_{1}\longrightarrow S_{X}\stackrel{\pi_{1}}{\longrightarrow} X\longrightarrow 0$
and $0\longrightarrow Y'\longrightarrow S_{Y}\stackrel{\pi_{2}}{\longrightarrow} Y\longrightarrow 0$ with $S_{X}, ~S_{Y}\in\mathrm{Add}S$
and $X_{1}, ~Y'\in \mathrm{Pres}^{n}_{\mathscr X}S$. We consider the following commutative diagram:
\begin{equation}\label{2.1}
\begin{split}
\xymatrix{
&0\ar[d]&0\ar[d]&0\ar[d]\\
0\ar[r]&X_{1}\ar[d]\ar[r]&Y_{1}\ar[d]\ar[r]&Z_{1}\ar[r]\ar[d]&0\\
0\ar[r]&S_{X}\ar[d]^{\pi_{1}}\ar[r]^(.4){\tiny \begin{pmatrix}1\\ 0\end{pmatrix}}
&S_{X}\bigoplus S_{Y}\ar[d]^{(i\pi_{1},\pi_{2})}\ar[r]^{~~~(0,1)}&S_{Y}\ar[r]\ar[d]^{\pi\pi_{2}}&0\\
0\ar[r]&X\ar[d]\ar[r]^{i}&Y\ar[d]\ar[r]^{\pi}&Z\ar[r]\ar[d]&0\\
&0&0&0
}
\end{split}\tag{$\natural$}
\end{equation}
It follows from the snake lemma that the first row in above diagram is $\mathscr X$-exact.
Since $S$ is an $n$-$\mathscr X$ quasi projective, the sequence
$0\longrightarrow Y'\longrightarrow S_{Y}\stackrel{\pi_{2}}{\longrightarrow} Y\longrightarrow 0$
is $S$-exact, and then the second column is $S$-exact. Thus $Y_{1}\in\mathrm{Pres}^{n}_{\mathscr X}S$ by Proposition \ref{char}.
In fact, all sequence in above are $\mathscr X$-exact.
Repeating the above process for the first row in above diagram, from the third column,
we can obtain that $Z\in\mathrm{Pres}^{n}_{\mathscr X}S$.

(2) If $X$ is in $\mathrm{Pres}^{n}_{\mathscr X}S$, the proof is similar to the Horse-shoe lemma.

If $Y$ is in $\mathrm{Pres}^{n}_{\mathscr X}S$, then there is an $\mathscr X$-exact sequence
$0\longrightarrow Y_{1}\longrightarrow S_{Y}\longrightarrow Y\longrightarrow 0$ with $S_{Y}\in\mathrm{Add}S$ and
$Y_{1}\in \mathrm{Pres}^{n}_{\mathscr X}S$. We consider the following communicative diagram.
$$\xymatrix{
&0\ar[d] & 0\ar[d]&\\
&Y_{1}\ar@^{=}[r]\ar[d] & Y_{1}\ar[d]\\
0\ar[r] &A\ar[r]\ar[d] &S_{Y}\ar[r]\ar[d] &Z\ar@^{=}[d]\ar[r] &0\\
0\ar[r]&X\ar[r]\ar[d] &Y\ar[d]\ar[r]&Z\ar[r] &0\\
&0&0&&}$$
The fact that the second row is both $S$-exact and $\mathscr X$-exact is not difficult to prove.
From the second row, we can get that $A\in\mathrm{Pres}^{n}_{\mathscr X}S$.
Note that the first column is $\mathscr X$-exact by Snake lemma.
It follows from the first column that $X\in\mathrm{Pres}^{n}_{\mathscr X}S$ by (1).

(3) We can obtain the commutative diagram ($\natural$) from the proof of (1), it is easy to see that $Z_{1}\in\mathrm{Pres}^{n}_{\mathscr X}S$ by (1),
and then the right column in diagram ($\natural$) is $S$-exact since $S$ is $(n+1)$-$\mathscr X$ quasi-projective. For any a morphism $f$: $S\longrightarrow Z$,
there is a morphism $g$: $S\longrightarrow S_{Y}$ such that $f=\pi\pi_{2} g$. It follows that $\mathrm{Hom}_{R}(S,\pi)$ is surjective.
i.e., the sequence $0\longrightarrow X\longrightarrow Y\longrightarrow Z\longrightarrow 0$ which is the third row in diagram ($\natural$) is $S$-exact.
\ \hfill $\Box$

\bg{Th}\label{xstarchar}
The following statements are equivalent:

$(1)$ S is an $n$-$\mathscr X$ star module;

$(2)$ If the exact sequence $0\longrightarrow U\longrightarrow V\longrightarrow W\longrightarrow 0$ is $\mathscr X$-exact
with $V, ~W\in \mathrm{Pres}^{n}_{\mathscr X}S$, then $U\in \mathrm{Pres}^{n}_{\mathscr X}S$ if and only if the sequence is S-exact.

$(3)$ If the sequence $0\longrightarrow U\longrightarrow S_{0}\longrightarrow V\longrightarrow 0$ is $\mathscr X$-exact
with $S_{0}\in\mathrm{Add}S$ and $V\in \mathrm{Pres}^{n}_{\mathscr X}S$, then $U\in \mathrm{Pres}^{n}_{\mathscr X}S$ if and only if the sequence
is S-exact;
\ed{Th}

\Pf. $(1)\Rightarrow (2)$ By (2) and (3) of the Proposition \ref{123}.

$(2)\Rightarrow (3)$ Obviously.

$(3)\Rightarrow (1)$ By the proposition \ref{char}, we only need to prove that $\mathrm{Pres}^{n}_{\mathscr X}S=\mathrm{Pres}^{n+1}_{\mathscr X}S$.
For any $M\in\mathrm{Pres}^{n}_{\mathscr X}S$, there is a $\mathscr X$-exact sequences
$0\longrightarrow M_{1}\stackrel{\alpha}{\longrightarrow} S_{M}\longrightarrow M\longrightarrow 0$
with $S_{M} \in\mathrm{Add}S$ and $M_{1} \in \mathrm{Pres}^{n-1}_{\mathscr X}S$. Then we have the following commutative diagram:
$$\xymatrix{
0\ar[r] &M_{1}\ar[r]^{\alpha}\ar[d]^{\gamma} &S_{M}\ar[r]^{g}\ar[d]^{f} &M\ar@^{=}[d]\ar[r] &0\\
0\ar[r]&M_{2}\ar[r] &S'_{M}\ar[r]^{\beta}&M\ar[r] &0
}$$
where $\beta$ is evaluation map. Indeed, the diagram above is a push-out of the morphism $\alpha$ and the morphism $\gamma$.
For any a morphism $\theta$: $X\longrightarrow M$ with $X\in \mathscr X$,
there is a morphism $\delta$: $X\longrightarrow S_{M}$ such that $\theta=g\delta$. i.e., $\theta=\beta f\delta$.
It follows that $\mathrm{Hom}_{R}(X,~\beta)$ is surjective. i.e., the second row is $\mathscr X$-exact.
Since the second row is $S$-exact by the structure, by the assumption, we have that $M_{2}\in\mathrm{Pres}^{n}_{\mathscr X}S$,
i.e., $M\in\mathrm{Pres}^{n+1}_{\mathscr X}S$. Thus we finish the proof.
\ \hfill $\Box$

\mskip\

Recall that a class of $R$-modules $\mathscr C$ is said to be closed under $\mathscr X$-extension, if for any $\mathscr X$-exact sequence
$0\longrightarrow X\longrightarrow Y\longrightarrow Z\longrightarrow 0$ with $X$, $Z\in \mathscr C$, then $Y\in \mathscr C$.

\bg{Pro}\label{xnchar}
The following statements are equivalent:

$(1)$ S is an $n$-$\mathscr X$ star module and $\mathrm{Pres}^{n}_{\mathscr X}S$ is closed under $\mathscr X$-extension;

$(2)$ $\mathrm{Pres}^{n}_{\mathscr X}S=\mathrm{Pres}^{n+1}_{\mathscr X}S\subseteq S^{\mathscr X\bot_{1}}$.

\ed{Pro}

\Pf. $(1)\Rightarrow (2)$ We only need to show that $\mathrm{Pres}^{n}_{\mathscr X}S\subseteq S^{\mathscr X\bot_{1}}$.
i.e., For any $M\in\mathrm{Pres}^{n}_{\mathscr X}S$, we only need to verify that the $\mathscr X$-exact sequence
$0\longrightarrow M\longrightarrow N\longrightarrow S\longrightarrow 0$
is split. Note that $N\in\mathrm{Pres}^{n}_{\mathscr X}S$ since $\mathrm{Pres}^{n}_{\mathscr X}S$ is closed under $\mathscr X$-extension.
It follows from Proposition \ref{123} (3) that the above sequence is $S$-exact, and then it is split.

$(2)\Rightarrow (1)$ Directly by definition, it's not hard to verify that $S$ is an $n$-$\mathscr X$ star module.
Now for any $\mathscr X$-exact $0\longrightarrow L\longrightarrow M\longrightarrow N\longrightarrow 0$ with $L,~N\in \mathrm{Pres}^{n}_{\mathscr X}S$.
Obviously, it is $S$-exact since $\mathrm{Pres}^{n}_{\mathscr X}S\subseteq S^{\mathscr X\bot_{1}}$.
By Proposition \ref{123} (2), $M\in\mathrm{Pres}^{n}_{\mathscr X}S$, i.e., $\mathrm{Pres}^{n}_{\mathscr X}S$ is closed under $\mathscr X$-extension.
\ \hfill $\Box$

\mskip\

Recall that a class of $R$-modules $\mathscr C$ is said to be closed under $n$-$\mathscr X$-images if for any $\mathscr X$-exact sequence
$C_{n}\longrightarrow \cdots \longrightarrow C_{2}\longrightarrow C_{1}\longrightarrow X\longrightarrow 0$ with $C_{i}\in\mathscr C$ for all $i$,
then $X\in\mathscr C$. Obviously, the class $\mathrm{Pres}^{1}_{\mathscr X}S$ is closed under $1$-$\mathscr X$-images. We do
not know that whether or not $\mathrm{Pres}^{n}_{\mathscr X}S$ is closed under $n$-$\mathscr X$-images in general. But if $S$ is a
$n$-$\mathscr X$ star module such that $\mathrm{Pres}^{n}_{\mathscr X}S$ is closed under $\mathscr X$-extensions, we have the following result.

\bg{Pro}\label{xnclosed}
Let S be an $n$-$\mathscr X$ star module such that $\mathrm{Pres}^{n}_{\mathscr X}S$ is closed under $\mathscr X$-extensions.
Then $\mathrm{Pres}^{k}_{\mathscr X}(\mathrm{Pres}^{n}_{\mathscr X}S)=\mathrm{Pres}^{k}_{\mathscr X}S$ for all $k\geq1$.
In particular, $\mathrm{Pres}^{n}_{\mathscr X}S$ is closed under $n$-$\mathscr X$-images.

\ed{Pro}

\Pf. We only need to verify that $\mathrm{Pres}^{k}_{\mathscr X}(\mathrm{Pres}^{n}_{\mathscr X}S)\subseteq\mathrm{Pres}^{k}_{\mathscr X}S$ for all $k\geq1$.
We proceed it by induction on $k$.

In case $k=1$, for any $A\in \mathrm{Pres}^{1}_{\mathscr X}(\mathrm{Pres}^{n}_{\mathscr X}S)$,
there is a $\mathscr X$-exact sequence
$0\longrightarrow A_{1}\longrightarrow B\longrightarrow A\longrightarrow 0$ with $B\in \mathrm{Pres}^{n}_{\mathscr X}S=\mathrm{Pres}^{n+1}_{\mathscr X}S$,
and then we have also a $\mathscr X$-exact sequence $0\longrightarrow B_{1}\longrightarrow S_{B}\longrightarrow B\longrightarrow 0$ with $S_{B} \in\mathrm{Add}S$
and $B_{1}\in\mathrm{Pres}^{n}_{\mathscr X}S$. We consider the following commutative diagram:
$$\xymatrix{
&0\ar[d] & 0\ar[d]&\\
&B_{1}\ar@^{=}[r]\ar[d] & B_{1}\ar[d]\\
0\ar[r] &C\ar[r]\ar[d] &S_{B}\ar[r]\ar[d] &A\ar@^{=}[d]\ar[r] &0\\
0\ar[r]&A_{1}\ar[r]\ar[d] &B\ar[d]\ar[r]&A\ar[r] &0\\
&0&0&&}$$
Note that the middle row is $\mathscr X$-exact. It follows that $A\in \mathrm{Pres}^{1}_{\mathscr X}S$.

Assume that the result holds for all $1\leq i\leq k$.
Next we will prove that the result is correct for $i=k+1$. Let $X\in\mathrm{Pres}^{k+1}_{\mathscr X}(\mathrm{Pres}^{n}_{\mathscr X}S)$, i.e.,
there is a $\mathscr X$-exact $Y_{k+1}\longrightarrow \cdots \longrightarrow Y_{2}\longrightarrow Y_{1}\stackrel{f}{\longrightarrow} X\longrightarrow 0$ with
$Y_{i}\in\mathrm{Pres}^{n}_{\mathscr X}S$ for all $i$. Set $\mathrm{Ker}f=X_{1}$, then $X_{1}\in\mathrm{Pres}^{k}_{\mathscr X}(\mathrm{Pres}^{n}_{\mathscr X}S)=\mathrm{Pres}^{k}_{\mathscr X}S$ by the induction assumption. Note that there is a $\mathscr X$-exact sequence
$0\longrightarrow Y_{1}'\longrightarrow S'\longrightarrow Y_{1}\longrightarrow 0$ with $S' \in\mathrm{Add}S$ and $Y_{1}'\in \mathrm{Pres}^{n}_{\mathscr X}S$
since $S$ is an $n$-$\mathscr X$ star module. We consider the following commutative diagram:
\begin{equation}\label{2.2}
\begin{split}
\xymatrix{
&0\ar[d] & 0\ar[d]&\\
&Y_{1}'\ar@^{=}[r]\ar[d] & Y_{1}'\ar[d]\\
0\ar[r] &Z\ar[r]\ar[d] &S'\ar[r]\ar[d] &X\ar@^{=}[d]\ar[r] &0\\
0\ar[r]&X_{1}\ar[r]\ar[d] &Y_{1}\ar[d]\ar[r]&X\ar[r] &0\\
&0&0&&}
\end{split}\tag{$\natural1$}
\end{equation}
It is easy to show that all rows and columns are $\mathscr X$-exact. We have a $\mathscr X$-exact sequence
$0\longrightarrow X_{1}'\longrightarrow S''\longrightarrow X_{1}\longrightarrow 0$ with $S'' \in\mathrm{Add}S$ and $X_{1}'\in \mathrm{Pres}^{k-1}_{\mathscr X}S$
since $X_{1}\in\mathrm{Pres}^{k}_{\mathscr X}S$. We consider the following commutative diagram:

\begin{equation}\label{2.3}
\begin{split}
\xymatrix{
&&0\ar[d] & 0\ar[d]&\\
&&X_{1}'\ar@^{=}[r]\ar[d] & X_{1}'\ar[d]\\
0\ar[r] &Y_{1}'\ar[r]\ar@^{=}[d] &M\ar[r]\ar[d] &S''\ar[d]\ar[r] &0\\
0\ar[r]&Y_{1}'\ar[r] &Z\ar[d]\ar[r]&X_{1}\ar[r]\ar[d] &0\\
&&0&0&}
\end{split}\tag{$\natural2$}
\end{equation}
It is easy to see that all rows and columns are $\mathscr X$-exact.
In the diagram (\ref{2.3}), from the middle row, we have that $M\in\mathrm{Pres}^{n}_{\mathscr X}S$ since
$\mathrm{Pres}^{n}_{\mathscr X}S$ is closed under $\mathscr X$-extensions.
Note that $X_{1}'\in \mathrm{Pres}^{k-1}_{\mathscr X}S=\mathrm{Pres}^{k-1}_{\mathscr X}(\mathrm{Pres}^{n}_{\mathscr X}S)$.
From the second column in (\ref{2.3}), we can know that $Z\in\mathrm{Pres}^{k}_{\mathscr X}(\mathrm{Pres}^{n}_{\mathscr X}S)=\mathrm{Pres}^{k}_{\mathscr X}S$.
It follows from the middle row in (\ref{2.2}) that $X\in\mathrm{Pres}^{k+1}_{\mathscr X}S$.
\ \hfill $\Box$

\section{Relative tilting modules}

In this section, we mainly give the relationship between $n$-$\mathscr X$ star modules and $n$-$\mathscr X$ tilting modules, see the theorem \ref{xtiltingstar}.
And we can obtain a new characterization of $n$-$\mathscr X$ tilting modules, see the theorem \ref{11}.

\bg{Def}$\mathrm{\cite[Definition ~3.1]{zpy}}$\label{xtilting}
An R-module T is said to be tilting with respect to $({\mathscr X}$, ${\mathscr Y})$ (or, call $n$-$\mathscr X$ tilting)
if it satisfying the following conditions:

$(1)$ ${\mathscr X}$-dim $T \leq n$;

$(2)$ $T$ is ${\mathscr X}$-self-orthogonal, i.e., $\mathrm{Ext}^{i}_{{\mathscr X}}(T,~T^{(I)})=0$, for each $i>0$ and all sets $I$;

$(3)$ there is a ${\mathscr X}$-exact sequence
$$\xymatrix{
0\ar[r]&X\ar[r]&T_{0}\ar[r]&\cdots\ar[r]&T_{n}\ar[r]&0
}$$
for any $X\in {\mathscr X}$, where $T_{i}\in \mathrm{Add}T$ for any $i$.
\ed{Def}

Next, we give the relationship between $n$-$\mathscr X$ star modules and $n$-$\mathscr X$ tilting modules.

\bg{Th}\label{xtiltingstar}
The following statements are equivalent:

$(1)$ $T$ is an $n$-$\mathscr X$ tilting module;

$(2)$ $\mathrm{Pres}^{n}_{\mathscr X}T=T^{\mathscr X\bot_{1\leq i\leq n}}$;

$(3)$ $T$ is an $n$-$\mathscr X$ star module and $\mathscr Y\subseteq \mathrm{Pres}^{n}_{\mathscr X}T$;

$(4)$ $\mathscr Y\subseteq \mathrm{Pres}^{n}_{\mathscr X}T=\mathrm{Pres}^{n+1}_{\mathscr X}T\subseteq T^{\mathscr X\bot_{1}}$.

\ed{Th}

\Pf. $(1)\Rightarrow(2)$ By \cite[Theorem~3.9]{zpy}.

$(2)\Rightarrow(3)$ By definition, we have that $\mathscr Y\subseteq T^{\mathscr X\bot_{1\leq i\leq n}}=\mathrm{Pres}^{n}_{\mathscr X}T$.
Next we only need to prove that $T$ is an $n$-$\mathscr X$ star module.
For any $\mathscr X$-exact sequence ($\ast$): $0\longrightarrow U\longrightarrow T' \longrightarrow V\longrightarrow 0$ with
$T'\in\mathrm{Add}T$ and $V\in\mathrm{Pres}^{n}_{\mathscr X}T$. It is easy to verify that $U\in T^{\mathscr X\bot_{2\leq i\leq n}}$ by the dimension shifting.
From the induced $\mathscr X$-exact sequence $0\longrightarrow \mathrm{Hom}_{R}(T,~U)\longrightarrow \mathrm{Hom}_{R}(T,T') \longrightarrow \mathrm{Hom}_{R}(T,~V)
\longrightarrow \mathrm{Ext}^{1}_{\mathscr X}(T,X)\longrightarrow 0$, we can obtain that the sequence ($\ast$) is $T$-exact if and only if
$U\in T^{\mathscr X\bot_{1}}$ if and only if $U\in T^{\mathscr X\bot_{1\leq i\leq n}}=\mathrm{Pres}^{n}_{\mathscr X}T$.
Consequently, $T$ is an $n$-$\mathscr X$ star module by Theorem \ref{xstarchar}

$(3)\Rightarrow(4)$ It is enough to prove that $\mathrm{Pres}^{n}_{\mathscr X}T\subseteq T^{\mathscr X\bot_{1}}$.
For any $U\in\mathrm{Pres}^{n}_{\mathscr X}T$, there exists $\mathscr X$-exact sequence ($\ast1$):
$0\longrightarrow U\longrightarrow Y \longrightarrow V\longrightarrow 0$ with
$Y\in\mathscr Y\subseteq \mathrm{Pres}^{n}_{\mathscr X}T$ by the assumption.
By the proposition \ref{123} (1), $V\in\mathrm{Pres}^{n}_{\mathscr X}T$. By the proposition \ref{123} (3), we know that the
sequence ($\ast1$) is $T$-exact. It follows from $\mathscr Y\subseteq T^{\mathscr X\bot_{1}}$ that $U\in T^{\mathscr X\bot_{1}}$.
i.e., $\mathrm{Pres}^{n}_{\mathscr X}T\subseteq T^{\mathscr X\bot_{1}}$.

$(4)\Rightarrow(1)$ We only need to prove that $\mathrm{Pres}^{n}_{\mathscr X}T=T^{\mathscr X\bot_{i\geq1}}$ by \cite[Theorem~3.9]{zpy}.
For any $U\in\mathrm{Pres}^{n}_{\mathscr X}T$, taking left $\mathscr Y$-coresolution of $U$, i.e., there exists $\mathscr X$-exact sequence
$0\longrightarrow U\longrightarrow Y_{0} \longrightarrow V\longrightarrow 0$ with $Y_{0}\in\mathscr Y\subseteq \mathrm{Pres}^{n}_{\mathscr X}T$.
By the proposition \ref{123} (1) and (3), $V\in\mathrm{Pres}^{n}_{\mathscr X}T$ and the sequence is $T$-exact.
It is easy to see that $U\in T^{\mathscr X\bot_{1,~2}}$. i.e.,$\mathrm{Pres}^{n}_{\mathscr X}T\subseteq T^{\mathscr X\bot_{1,2}}$.
Repeat the above process for $V$, we can prove that $\mathrm{Pres}^{n}_{\mathscr X}T\subseteq T^{\mathscr X\bot_{i\geq1}}$.

On the other hand, for any $M\in T^{\mathscr X\bot_{\geq1}}$, taking left $\mathscr Y$-coresolution of $M$, i.e., there is a long $\mathscr X$-exact sequence
$0\longrightarrow M\stackrel{f_{0}}{\longrightarrow} Y_{0}\stackrel{f_{1}}{\longrightarrow} Y_{1}\longrightarrow\cdots
\stackrel{f_{n}}{\longrightarrow} Y_{n}\longrightarrow Y_{n+1}$ with
$Y_{j}\in\mathscr Y$ for all $j$. Set $\mathrm{Ker}f_{j}=M_{j-1}$, $1\leq j\leq n$, where
$M_{0}=M$. By the dimension shifting, $M_{j}\in T^{\mathscr X\bot_{\geq1}}$ for all $j$.
By the proposition \ref{xnchar}, $\mathrm{Pres}^{n}_{\mathscr X}T$ is closed under $\mathscr X$-extension.
Since $\mathscr Y\subseteq \mathrm{Pres}^{n}_{\mathscr X}T$,
$M_{n}\in\mathrm{Pres}^{n}_{\mathscr X}T$. By the proposition \ref{123} (2), $M\in\mathrm{Pres}^{n}_{\mathscr X}T$.
i.e., $T^{\mathscr X\bot_{i\geq1}}\subseteq\mathrm{Pres}^{n}_{\mathscr X}T$.
\ \hfill $\Box$

\mskip\

Let $T$ be an $R$-module. Denote $\mathcal{C}_{n}^{\mathscr XT}$ the subcategory there is an $\mathscr X$-exact sequence
$T_{n}\longrightarrow T_{n-1}\cdots\longrightarrow T_{1}\longrightarrow T_{0} \longrightarrow U\longrightarrow0$
with $T_{i}\in \mathrm{Add}T$, and the exact sequence is $T$-exact.
In order to give a new characterization of $n$-$\mathscr X$ tilting modules,
we need the following two lemmas.

\bg{Lem}\label{cngt}
Let $0\longrightarrow U\stackrel{i}{\longrightarrow} V\stackrel{\pi}{\longrightarrow} W\longrightarrow 0$ be an $\mathscr X$-exact sequence.

$(1)$ If U and V are in $\mathcal{C}_{n}^{\mathscr XT}$, then $W\in \mathrm{Pres}^{n+1}_{\mathscr X}T$. Moveover, if the sequence is
T-exact, then $W\in \mathcal{C}_{n}^{\mathscr XT}$;

$(2)$ If V and W are in $\mathcal{C}_{n}^{\mathscr XT}$ and the sequence is T-exact, then $U\in \mathcal{C}_{n-1}^{\mathscr XT}$.
\ed{Lem}

\Pf. $(1)$ We show it by induction on $n$. For $n=0$, the result is obvious.

Now, we assume that the conclusion holds for $n=k-1$. Let $U$, $V \in \mathcal{C}_{k}^{\mathscr XT}$, we have two $\mathscr X$-exact sequences
$0\longrightarrow U'\longrightarrow T_{U}\stackrel{f}{\longrightarrow} U\longrightarrow 0$
and $0\longrightarrow V'\longrightarrow T_{V}\stackrel{g}{\longrightarrow} V\longrightarrow 0$ with $T_{U}, ~T_{V}\in\mathrm{Add}T$
and $U', ~V'\in \mathcal{C}_{k-1}^{\mathscr XT}$ such that both sequences are $T$-exact.
We consider the following commutative diagram:
$$\xymatrix{
&0\ar[d]&0\ar[d]&0\ar[d]\\
0\ar[r]&U'\ar[d]\ar[r]&V''\ar[d]\ar[r]&W_{1}\ar[r]\ar[d]&0\\
0\ar[r]&T_{U}\ar[d]\ar[r]^(.4){\tiny \begin{pmatrix}1\\ 0\end{pmatrix}}
&T_{U}\bigoplus T_{V}\ar[d]^{(if,g)}\ar[r]^{~~~(0,1)}&T_{V}\ar[r]\ar[d]^{\pi g}&0\\
0\ar[r]&U\ar[d]\ar[r]^{i}&V\ar[d]\ar[r]^{\pi}&W\ar[r]\ar[d]&0\\
&0&0&0
}$$
It is easy to see that all sequences in above diagram are $\mathscr X$-exact. From the first row, we can obtain a new $\mathscr X$-exact sequence
$0\longrightarrow U'\bigoplus T_{V}\longrightarrow V''\bigoplus T_{V}\longrightarrow W_{1}\longrightarrow 0$. By Lemma 2.3 in \cite{WEI2},
from the two $\mathscr X$-exact sequences $0\longrightarrow V'\longrightarrow T_{V}\stackrel{g}{\longrightarrow} V\longrightarrow 0$
and $0\longrightarrow V''\longrightarrow T_{U}\bigoplus T_{V}\longrightarrow V\longrightarrow 0$, we have that
$V''\bigoplus T_{V}\cong V'\bigoplus T_{U}\bigoplus T_{V}\in \mathcal{C}_{k-1}^{\mathscr XT}$. It follows from the induction assumption
that $W_{1}\in\mathrm{Pres}^{k}_{\mathscr X}T$ for the first row in above diagram. From the right column, we have that $W\in\mathrm{Pres}^{k+1}_{\mathscr X}T$.
Moveover, if the sequence is $T$-exact, it is not difficult to verify that all sequences in diagram above
is also $T$-exact, then the next proof is easy.

$(2)$ Since both $V$ and $W$ are in $\mathcal{C}_{n}^{\mathscr XT}$, there are two $\mathscr X$-exact sequences
$0\longrightarrow V'\longrightarrow T_{V}\longrightarrow V\longrightarrow 0$
and $0\longrightarrow W'\longrightarrow T_{W}\longrightarrow W\longrightarrow 0$ with $T_{V}, ~T_{W}\in\mathrm{Add}T$
and $V', ~W'\in \mathcal{C}_{n-1}^{\mathscr XT}$ such that the two sequences are $T$-exact.
We consider the following commutative diagram:
$$\xymatrix{
&0\ar[d] & 0\ar[d]&\\
&V'\ar@^{=}[r]\ar[d] & V'\ar[d]\\
0\ar[r] &W''\ar[r]\ar[d] &T_{V}\ar[r]\ar[d] &W\ar@^{=}[d]\ar[r] &0\\
0\ar[r]&U\ar[r]\ar[d] &V\ar[d]\ar[r]&W\ar[r] &0\\
&0&0&&}$$
It is easy to verify that all sequences in diagram above are $\mathscr X$-exact and $T$-exact.
From the left column, we can obtain a new $\mathscr X$-exact sequence ($\ast$2):
$0\longrightarrow V'\bigoplus T_{W}\longrightarrow W''\bigoplus T_{W}\longrightarrow U\longrightarrow 0$
with $V'\bigoplus T_{W}\in \mathcal{C}_{n-1}^{\mathscr XT}$. Clearly, the new sequence is $T$-exact.
Consider the second column and the sequence $0\longrightarrow W'\longrightarrow T_{W}\longrightarrow W\longrightarrow 0$,
by Lemma 2.3 in \cite{WEI2}, we have that $W''\bigoplus T_{W}\cong W'\bigoplus T_{V}\in \mathcal{C}_{n-1}^{\mathscr XT}$.
From the sequence ($\ast$2), we have that $U\in \mathcal{C}_{n-1}^{\mathscr XT}$ by (1).
\ \hfill $\Box$

\bg{Lem}\label{add}
If $T$ is an $n$-$\mathscr X$ tilting module, then $\mathcal{C}_{n}^{\mathscr XT}=\mathrm{Pres}^{n}_{\mathscr X}T$.
\ed{Lem}

\Pf. Clearly, $\mathcal{C}_{n}^{\mathscr XT}\subseteq \mathrm{Pres}^{n+1}_{\mathscr X}T=\mathrm{Pres}^{n}_{\mathscr X}T$.
On the other hand, for any $M\in\mathrm{Pres}^{n+1}_{\mathscr X}T$, there is an infinite $\mathscr X$-exact sequence
\begin{equation}\label{3.1}
\begin{split}
\xymatrix{\cdots\ar[r]&T_{n}\ar[r]^{f_{n}}&\cdots \ar[r]&T_{1}\ar[r]^{f_{1}}&T_{0}\ar[r]^{f_{0}}&M\ar[r]&0}
\end{split}\tag{$\natural3$}
\end{equation}
with $T_{i}\in \mathrm{Add}T$ for all $i$. Set $\mathrm{Im}f_{i}=M_{i}\in\mathrm{Pres}^{n+1}_{\mathscr X}T$.
Consider the $\mathscr X$-exact sequence ($\ast$3): $\xymatrix{0\ar[r]&M_{1}\ar[r]&T_{0}\ar[r]&M\ar[r]&0}$.
By the theorem \ref{xtiltingstar}, $T$ is $n$-$\mathscr X$ quasi-projective. Then the sequence ($\ast$3) is
$T$-exact by the proposition \ref{char}. From the sequence ($\natural$3), we can know that $M\in\mathcal{C}_{n}^{\mathscr XT}$.
So $\mathcal{C}_{n}^{\mathscr XT}= \mathrm{Pres}^{n+1}_{\mathscr X}T=\mathrm{Pres}^{n}_{\mathscr X}T$.
\ \hfill $\Box$

\mskip\

At the end of this paper, we can give a new characterization of $n$-$\mathscr X$ tilting modules.

\bg{Th}\label{11}
Let T be an R-module. Then T is an $n$-$\mathscr X$ tilting module if and only if

$(1)$ $\mathrm{Hom}(T,-)$ preserve exactness in $\mathcal{C}_{n}^{\mathscr XT}$,

$(2)$ $\mathscr Y\subseteq \mathcal{C}_{n}^{\mathscr XT}$,

$(3)$ $\mathcal{C}_{n}^{\mathscr XT}$ is closed under $n$-$\mathscr X$-images.
\ed{Th}

\Pf. ($\Rightarrow$) By Theorem \ref{xtiltingstar}, Lemma \ref{add}, Proposition \ref{xnchar} and Proposition \ref{xnclosed}.

($\Leftarrow$) Note that $\mathcal{C}_{n}^{\mathscr XT}\subseteq \mathrm{Pres}^{n+1}_{\mathscr X}T\subseteq\mathrm{Pres}^{n}_{\mathscr X}T$. Since
$\mathcal{C}_{n}^{\mathscr XT}$ be closed under $n$-$\mathscr X$-images and $\mathrm{Add}T\in\mathcal{C}_{n}^{\mathscr XT}$,
we have that $\mathrm{Pres}^{n}_{\mathscr X}T\subseteq\mathcal{C}_{n}^{\mathscr XT}$.
i.e., $\mathcal{C}_{n}^{\mathscr XT}= \mathrm{Pres}^{n+1}_{\mathscr X}T=\mathrm{Pres}^{n}_{\mathscr X}T$.
For any $M\in \mathcal{C}_{n}^{\mathscr XT}$, there is an $\mathscr X$-exact sequence ($\ast4$):
$0\longrightarrow M\longrightarrow Y\longrightarrow M_{1}\longrightarrow0$ with $Y\in \mathscr Y$.
By the lemma \ref{cngt} (1), $M_{1}\in \mathrm{Pres}^{n+1}_{\mathscr X}T=\mathcal{C}_{n}^{\mathscr XT}$.
Consequently, the sequence ($\ast4$) is $T$-exact by the assumption, and then $M\in T^{\mathscr X\bot_{1}}$, i.e.,
$\mathcal{C}_{n}^{\mathscr XT}\subseteq T^{\mathscr X\bot_{1}}$. So $T$ is an $n$-$\mathscr X$ tilting module by Theorem \ref{xtiltingstar} (4).
\ \hfill $\Box$

%

{\small

}

\end{document}